\input graphicx
\input amstex
\input amssym.def
\UseAMSsymbols
\documentstyle{amsppt}
\def\EZ{even-zeroed}
\def\nagyEZ{Even-zeroed}
\nopagenumbers
\topmatter
\title A simple combinatorial proof of Shapiro's Catalan convolution
\endtitle
\author G\'abor V. Nagy\endauthor
\affil Bolyai Institute, University of Szeged, Hungary\endaffil
\abstract
Shapiro proved an elegant convolution formula involving Catalan numbers of even index.
This paper gives a simple combinatorial proof of his formula.
In addition, we show that it is equivalent with the alternating convolution formula of central
binomial coefficients.
\endabstract
\endtopmatter
\document
\specialhead1. Introduction\endspecialhead
In this paper
\footnote""{
This research was partially supported by the Hungarian Scientic
Research Fund, Grant No. K
75517 and by the TAMOP-4.2.2/B-10/1-2010-0012 and
TAMOP-4.2.1/B-09/1/KONV-2010-0005
programs of the Hungarian National Development Agency.\hfil
}
$C_n$ denotes the $n$th Catalan-number and $B_n$ denotes the $n$th central
binomial coefficient, i.e. $C_n=\frac1{n+1}\binom{2n}n$ and $B_n=\binom{2n}n$.
Unless otherwise stated, all indices (i, j, k, and so on) are nonnegative integers in our formulas.\par
In 2002, L. Shapiro found the following elegant identity [3; p. 123]:
\proclaim{Theorem 1}
$$\sum_{i+j=n}C_{2i}C_{2j}=4^nC_n\eqno(1)$$
\endproclaim
This can be easily proved using generating functions, but according to Stanley [5; p. 46], no simple combinatorial proof has been known
(see [1] for a recent complicated combinatorial proof). In Section 3 we will give a simple combinatorial proof
of the the following equivalent version of Theorem 1:
\proclaim{Theorem 2}
$$\sum_{i+j=n}C_{2i}B_{2j}=4^nB_n\eqno(2)$$
\endproclaim
\noindent(1) and (2) are equivalent, because
$$\aligned\sum_{i+j=n}C_{2i}B_{2j}&=
\frac12\left(\sum_{i+j=n}C_{2i}(2j+1)C_{2j}+\sum_{i+j=n}(2i+1)C_{2i}C_{2j}\right)=\\
&=\frac12\sum_{i+j=n}(2n+2)C_{2i}C_{2j}=(n+1)\sum_{i+j=n}C_{2i}C_{2j}.\endaligned$$\par
The key observation of this paper is a non-standard interpretation of $C_{2n}$, that is discussed in the next section.
Using that, we can give a new combinatorial meaning of the left-hand sides of (1) and (2).
In Section 4, we show bijectively that (2) is equivalent with the alternating convolution formula of central
binomial coefficients, which has a nice combinatorial interpretation, due to Spivey [4].
\specialhead2. \nagyEZ\ balanced paths and $C_{2n}$\endspecialhead
A {\it path of length $l$\/} is an $l$-element sequence of up-steps ($\nearrow$) and down-steps ($\searrow$).
A {\it balanced $n$-path\/} is such a path of length $2n$ that has $n$ up-steps and $n$ down-steps.
The number of balanced $n$-paths is clearly $B_n$. We denote by $\Cal B_n$ the set of balanced $n$-paths.
We visualize paths in the usual way, they start from the origin,
$\nearrow$ is a step $(1,1)$, $\searrow$ is a step $(1,-1)$, see the figures below.
An {\it $n$-Dyck-path\/} (of length $2n$)
is a balanced $n$-path such that it never falls below the $x$-axis. 
It is well-known that the number of $n$-Dyck-paths is $C_n$. We denote by $\Cal C_n$ the set of $n$-Dyck-paths.
A {\it signed $n$-Dyck-path\/} is an element of the set $\Cal S_n:=\{+,-\}\times\Cal C_n$. The set of 
signed Dyck-paths is denoted by $\Cal S:=\cup_{i=0}^\infty\Cal S_i$.
The number of up-steps (or down-steps) in a signed or unsigned balanced path $P$
is called the parameter of $P$ and it is denoted by $\text{par}(P)$.
A special type of paths plays crucial role in our proofs, so we introduce a new terminology for them:
We call a (balanced or non-balanced) path {\it\EZ\/}, if its $x$-intercepts are all divisible by $4$.\par
If $n\ge1$, then every balanced $n$-path can be decomposed uniquely into a sequence of signed Dyck-paths in a very natural way (see Figure 1):
The $x$-axis cuts the balanced path into nonempty subpaths so that every subpath is either a Dyck-path that never touches the $x$-axis (apart from its starting point and end point),
or the reflection of such a Dyck-path across the $x$-axis. Every subpath is uniquely characterized by a signed Dyck-path that we
get after removing the first and last steps ($+$: standard Dyck-path, $-$: reflected Dyck-path), and we can list these signed Dyck-paths (from left to right) in a sequence.
It is very easy to see that we defined a $\chi$ bijection between $\Cal B_n$ and $\text{SEQ}(n)$, where
$$\text{SEQ}(n)=\left\{(P_1,\dots,P_k): k\in\Bbb Z^+, P_i\in\Cal S\text{ for all } i,\text{ and }\sum_{i=1}^k(\text{par}(P_i)+1)=n\right\}.$$
Now we are ready to prove the key lemma of this paper:
\proclaim{Lemma 3}
$C_{2n}$ counts the number of \EZ\ balanced $2n$-paths.
\endproclaim
\demo{Proof}
The statement is true for $n=0$. Now let us assume that $n\ge1$.\par
Clearly, a balanced $2n$-path $P$ is \EZ, if and only if
all the signed Dyck-paths has odd parameter in its sequence $\chi(P)$. So if we denote by $\tilde\Cal B_{2n}$ the set of \EZ\ $2n$-paths,
the restriction of $\chi$ to $\tilde\Cal B_{2n}$ gives a $\phi$ bijection between $\tilde\Cal B_{2n}$ and $\widetilde{\text{SEQ}}(2n)$, where
$$\widetilde{\text{SEQ}}(2n)=\left\{(P_1,\dots,P_k): P_i\in\Cal S,\text{ par}(P_i)\text{ is odd } \forall i;\ \sum_{i=1}^k(\text{par}(P_i)+1)=2n\right\}.$$\par
Now we define a $\psi$ bijection between $\Cal C_{2n}$ and $\widetilde{\text{SEQ}}(2n)$ and so we give a bijective proof of the lemma
($\phi^{-1}\circ\psi$ is a bijection between $\Cal C_{2n}$ and $\tilde\Cal B_{2n}$).
Consider an arbitrary $2n$-Dyck-path $D$. It is well-known from a standard proof of the Catalan-recursion that $D$ can be uniquely written as $\nearrow$$L$$\searrow$$R$,
i.e. $D$ can be decomposed into an $(L,R)$ ordered pair, where $L$ and $R$ are Dyck-paths, whose parameters sum to $2n-1$.
Either $\text{par}(L)$ or $\text{par}(R)$ is odd. If $\text{par}(L)$ is odd, then we define the first element of $\psi(D)$ to be $-L$,
and we recursively repeat the process for $D':=R$ ($\text{par}(R)$ is even) to get the other elements of $\psi(D)$.
If $\text{par}(R)$ is odd, then we define the first element of the $\psi(D)$ to be $+R$,
and recursively repeat the process for $D':=L$. ($-$ means ,,left'', $+$ means ,,right'' here.) The process terminates when $D'$ is the empty $0$-Dyck-path.
It is easy to check that the obtained $\psi(D)$ is in $\widetilde{\text{SEQ}}(2n)$. See Figure 2 for a visualization.\par
One can easily compute $\psi^{-1}(S)$ for an arbitrary $S\in\widetilde{\text{SEQ}}(2n)$, so we indeed defined a bijection.\qed
\enddemo
\centerline{\includegraphics[width=331pt]{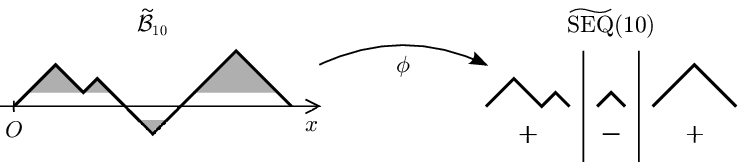}}
\centerline{\eightbf Figure 1: Illustration of $\chi$ and $\phi$}\smallskip
\centerline{\includegraphics[width=331pt]{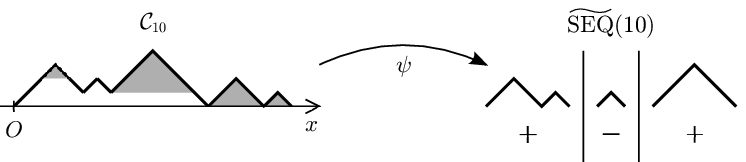}}
\centerline{\eightbf Figure 2: Illustration of $\psi$}\par
\demo{Remarks} Roughly speaking, our bijection $\Cal C_{2n}\to\tilde\Cal B_{2n}$ converts the ,,left-right symmetry'' of $\Cal C_{2n}$
into the ,,up-down symmetry'' of $\tilde\Cal B_{2n}$.\par
When defining $\psi$, it might be slightly more natural to work with full binary trees that is an other representation of Catalan-numbers
(what we do here is to decompose even-parameter full binary trees into odd-parameter subtrees -- the details are left to reader). But
then $\phi^{-1}\circ\psi$ would become slightly less intuitive, since we would need an extra conversion between full binary trees and Dyck-paths.\par
If we already know or conjecture that $|\tilde\Cal B_{2n}|=C_{2n}$, we can find a quicker (but recursive) argument for this. Namely,
using the notations $X_n:=|\tilde\Cal B_{2n}|$ and $Y_n:=C_{2n}$, one can quickly figure out that both $(X_n)_{n=0}^\infty$ and $(Y_n)_{n=0}^\infty$
satisfy the following recursion: $Z_0=1$,\quad$Z_n=2\sum_{k=1}^nC_{2k-1}Z_{n-k}$ (if $n\ge1$).
\enddemo
As an application, we prove a lemma, from which one can see that both sides of (2) satisfy the following recursion:
$X_0=1,\quad\sum_{i+j=n}X_iX_j=16^n$.
\proclaim{Lemma 4} $$2\cdot\sum_{i+j+k=n}C_{2i}C_{2j}B_{2k}=B_{2n+1}$$\endproclaim
\demo{Proof}
Both sides count the number of balanced $(2n+1)$-paths. This is obvious for the right-hand side.
In the left-hand side, we group the balanced $(2n+1)$-paths by the position (i) and length (j)
of the leftmost such signed Dyck-path segment (cut by the $x$-axis) whose parameter is odd.
Such a segment must exist, since the sum of the parameters is $2n+1$. If the starting point of that segment is $4i$ and its parameter is $2j+1$,
then by Lemma 3, there are $C_{2i}$ \EZ\ balanced paths from the origin to $4i$,
there are $2C_{2j}$ possible choices for the segment in question, and there are $B_{2k}$ possible endings for the rest of the path ($k=n-i-j$).\qed
\enddemo
\specialhead 3. The proof of Theorem 2\endspecialhead
The following lemma is well-known, and it has several combinatorial proofs~[2].
\proclaim{Lemma 5}$B_n$ counts the number of such paths of length $2n$ that never return to the $x$-axis after the first step.\endproclaim
With the help of Lemma 3, we can give an interesting combinatorial interpretation of the left-hand sides of (1) and (2).
\proclaim{Lemma 6}a) $\sum_{i+j=n}C_{2i}B_{2j}$ is the number of \EZ\ paths of length $4n$.\par
\noindent b) $\sum_{i+j=n}C_{2i}C_{2j}$ is the number of \EZ\ paths from the origin to $(4n+1,1)$.
\endproclaim
\demo{Proof}a) By Lemmas 3 and 5, there are $C_{2i}B_{2(n-i)}$ such \EZ\ paths of length $4n$ whose rightmost $x$-intercept is $4i$.\par
\noindent b) There are $C_{2i}C_{2(n-i)}$ such \EZ\ paths from the origin to $(4n+1,1)$ whose rightmost $x$-intercept is $4i$ (followed by an up-step).
The first factor comes from Lemma 3 and the second one comes from the standard interpretion of $C_{2(n-i)}$ by Dyck-paths.\qed
\enddemo
\centerline{\includegraphics[width=135pt]{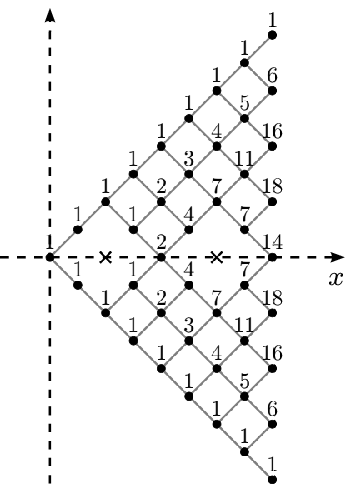}}
\centerline{\eightbf Figure 3: The number of even-zeroed paths}\smallskip
In Figure 3 the label of a node shows the number of \EZ\ paths from the origin to that node. These labels can be calculated recursively, since every label
is the sum of its left neighbors. We already know that the label of $(4n,0)$ is $C_{2n}$, the label of $(4n+1,\pm1)$ is $L_n:=\sum_{i+j=n}C_{2i}C_{2j}$
and the sums of the labels in the $4n$th row is $S_n:=\sum_{i+j=n}C_{2i}B_{2j}$. In order to prove Theorem 2, we only have to show that $S_n=4^nB_n$.
The key observation is that $S_{n+1}$ can be calculated from $S_n$ and $L_n$ easily, but we know from Section 1, that $L_n=\frac1{n+1}S_n$, so in fact
$S_{n+1}$ can be calculated from $S_n$ easily. This calculation is done in the next lemma, which implies Theorem 2.
\proclaim{Lemma 7} The number of \EZ\ paths of length $4n$ is $4^nB_n$.
\endproclaim
\demo{Proof}
Let $\Cal P_n$ denote the set of \EZ\ paths of length $4n$, and set $S_n:=|\Cal P_n|$. By induction on $n$, we prove that $S_n=4^nB_n$.
This is obviously true if $n=0$.\par
Let us assume that $S_n=4^nB_n$ holds. Clearly, 
every path of $\Cal P_{n+1}$ is an extension of a path of $\Cal P_n$
by $4$ steps. For each path of $\Cal P_n$ there are $16$ possible extensions.
But some of the $16S_n$ extensions are not in $\Cal P_{n+1}$. These ,,wrong'' extenstions are exactly the \EZ\ paths from the origin to $(4n+1,1)$
followed by a down-step and two arbitrary steps, and the reflections of these paths across the $x$-axis.
By Lemma 6.b, the number of these wrong extensions is $8\sum_{i+j=n}C_{2i}C_{2j}$, that equals to $\frac8{n+1}\sum_{i+j=n}C_{2i}B_{2j}=\frac8{n+1}S_n$,
as seen in Section 1 and Lemma 6.a.
By the induction hypothesis, $S_n=4^nB_n$, thus $S_{n+1}=16\cdot4^nB_n-\frac8{n+1}4^nB_n$.
A quick calculation shows that $S_{n+1}=4^{n+1}B_{n+1}$.\qed
\enddemo
\specialhead 4. Alternating convolution of the central binomial coefficients\endspecialhead
The following theorem has a nice combinatorial proof, due to Spivey [4]:
\proclaim{Theorem 8}
$$\sum_{i+j=n}B_{2i}B_{2j}-\sum\Sb i+j=n\\j\ge1\endSb B_{2i+1}B_{2j-1}=4^nB_n.$$
\endproclaim
By proving the next theorem bijectively, we will see that Theorem 8 is equivalent with Theorem 2, so any combinatorial proof of Theorem 8
yields a combinatorial proof of Theroem 2. Conversely, our proof in the previous section can be interpreted as a new proof of Theorem 8.
\proclaim{Theorem 9}
$$\sum_{i+j=n}B_{2i}B_{2j}-\sum\Sb i+j=n\\j\ge1\endSb B_{2i+1}B_{2j-1}=\sum_{i+j=n}C_{2i}B_{2j}.$$
\endproclaim
\demo{Proof}Using Lemma 3, we will prove the following equivalent form:
$$\sum_{i+j=n}(B_{2i}-C_{2i})B_{2j}=\sum\Sb i+j=n\\j\ge1\endSb B_{2i+1}B_{2j-1}.\eqno(3)$$\par
The right-hand side counts the number of $(O_1,O_2)$ pairs, where $O_1$ and $O_2$ are balanced paths with odd parameters,
and $\text{par}(O_1)+\text{par}(O_2)=2n$. Let $\Cal O$ be the set of these pairs.
By Lemma 3, the left-hand side counts the number of $(E_1,E_2)$ pairs, where
$E_1$ and $E_2$ are balanced paths with even parameters, $E_1$ has an $x$-intercept of the form $4t+2$ ($t\in\Bbb N_0$),
and $\text{par}(E_1)+\text{par}(E_2)=2n$. Let $\Cal E$ be the set of these pairs.\par
We will give a bijection between $\Cal E$ and $\Cal O$, which means that $|\Cal E|=|\Cal O|$, as stated.
Pick an arbitrary element $(E_1, E_2)$ of $\Cal E$.
Let $L$ be such subpath of $E_1$ that is identical with $E_1$ from the origin to its leftmost $x$-intercept of the form $4t+2$, and let $R$ be the rest of $E_1$.
Then the image of $(E_1, E_2)$ is defined as $(LE_2,R)$, where $LE_2$ is the concatenation of $L$ and $E_2$ in this order.
It is easy to see that this mapping is bijective.\qed
\enddemo
If we write $C_{2i}=\binom{4i}{2i}-\binom{4i}{2i-1}$ in (3), we get the following identity:
\proclaim{Corollary 10}
$$\sum_{i=1}^{n}\binom{4i}{2i-1}\binom{4n-4i}{2n-2i}=\sum_{i=0}^{n-1}\binom{4i+2}{2i+1}\binom{4n-4i-2}{2n-2i-1}$$
\endproclaim
\Refs
\item{[1]} G. E. Andrews, On Shapiro's Catalan convolution, Advances in Applied Mathematics, Volume 46, Issues 1–4, January 2011, Pages 15-24.
\item{[2]} \"O. E\u gecio\u glu and A. King, Random walks and Catalan factorization, Congressus Numerantium, 138 (1999), Pages 129-140.
\item{[3]} T. Koshy, Catalan Numbers with Applications, Oxford University Press, New York, 2009.
\item{[4]} M. Spivey, Combinatorial interpretation of the alternating convolution of the central binomial coefficients,
http://mikespivey.wordpress.com/2012/03/16/altconvcentralbinom/, March 16, 2012.
\item{[5]} R. P. Stanley, Bijective proof problems, http://www-math.mit.edu/\~{}rstan/bij.pdf, August 18, 2009.
\endRefs
\enddocument